\documentclass[12pt, a4paper,reqno]{amsart}
\usepackage{amsthm}
\usepackage{mathtools}
{\theoremstyle{definition}
\newtheorem{dfn}{Definition}[section]}
\newtheorem{prop}[dfn]{Proposition}

\newtheorem{thm}[dfn]{Theorem}
\theoremstyle{definition}
\newtheorem{rem}[dfn]{Remark}
\newtheorem{lem}[dfn]{Lemma}

{\theoremstyle{definition}
}

\newtheorem{prob}[dfn]{Problem}
\usepackage[top=30truemm,bottom=30truemm,left=35truemm,right=35truemm]{geometry}
\usepackage{amsmath,amssymb,amscd,bm,graphicx,tikz-cd,upgreek,dsfont,amsfonts,tikz,enumitem}
\usetikzlibrary{cd}
\usepackage[all]{xy}
\usepackage[colorlinks=true,
    citecolor=green(pigment),
    linkcolor=alizarin,
    urlcolor=denim]{hyperref}
    
\definecolor{alizarin}{rgb}{0.82, 0.1, 0.26}
\definecolor{azure(colorwheel)}{rgb}{0.0, 0.5, 1.0}
\definecolor{blue(pigment)}{rgb}{0.2, 0.2, 0.6}
\definecolor{denim}{rgb}{0.08, 0.38, 0.74}
\definecolor{mint}{rgb}{0.24, 0.71, 0.54}
\definecolor{parisgreen}{rgb}{0.31, 0.78, 0.47}
\definecolor{persiangreen}{rgb}{0.0, 0.65, 0.58}
\definecolor{seagreen}{rgb}{0.18, 0.55, 0.34}
\definecolor{shamrockgreen}{rgb}{0.0, 0.62, 0.38}
\definecolor{green(pigment)}{rgb}{0.0, 0.65, 0.31}

\usepackage{mathrsfs}

\usepackage{setspace}
\newcommand{\marginparstretch}{0.6}
\let\oldmarginpar\marginpar
\renewcommand\marginpar[1]{\-\oldmarginpar[\framebox{\setstretch{\marginparstretch}\begin{minipage}{\marginparwidth}{\raggedleft\tiny #1}\end{minipage}}]{\framebox{\setstretch{\marginparstretch}\begin{minipage}{\marginparwidth}{\raggedright\tiny #1}\end{minipage}}}}

\usepackage{scalerel,stackengine}
\stackMath
\newcommand\reallywidecheck[1]{%
\savestack{\tmpbox}{\stretchto{%
   \scaleto{%
    \scalerel*[\widthof{\ensuremath{#1}}]{\kern-.6pt\bigwedge\kern-.6pt}%
    {\rule[-\textheight/2]{1ex}{\textheight}}
   }{\textheight}%
}{0.5ex}}%
\stackon[1pt]{#1}{\scalebox{-1}{\tmpbox}}%
}
\def\ang<#1>{\langle #1 \rangle}
\def\bigang<#1>{\left\langle #1 \right\rangle}

\numberwithin{equation}{section}

\newcommand{\Spec}{\operatorname{Spec}}

\newcommand{\Hom}{\operatorname{Hom}}

\newcommand{\Pic}{\operatorname{Pic}}

\newcommand{\Db}{\operatorname{D}^{\mathrm{b}}}

\newcommand{\cI}{\mathcal{I}}

\newcommand{\cL}{\mathcal{L}}

\newcommand{\cO}{\mathcal{O}}

\newcommand{\cQ}{\mathcal{Q}}

\newcommand{\cX}{\mathcal{X}}
\newcommand{\cY}{\mathcal{Y}}

\newcommand{\bC}{\mathbb{C}}

\newcommand{\bP}{\mathbb{P}}

\newcommand{\bZ}{\mathbb{Z}}

\newcommand{\m}{\mathfrak{m}}

\newcommand{\simto}{\xrightarrow{\sim}}

\tikzset{
        DB/.style={circle,draw=black,circle,fill=white,inner sep=0pt, minimum size=4pt},
        DW/.style={circle,draw=black,fill=black,inner sep=0pt, minimum size=4pt},
        cvertex/.style={circle,draw=black,fill=white,inner sep=1pt,outer sep=3pt},
        vertex/.style={circle,fill=black,inner sep=1pt,outer sep=3pt},
        star/.style={circle,fill=yellow,inner sep=0.75pt,outer sep=0.75pt},
        tvertex/.style={inner sep=1pt,font=\scriptsize},
pvertex/.style={circle,inner sep=1pt,outer sep=2pt,font=\scriptsize},
        gap/.style={inner sep=0.5pt,fill=white}
}

\makeatletter
\newcommand*{\defeq}{\mathrel{\rlap{%
                     \raisebox{0.3ex}{$\m@th\cdot$}}%
                     \raisebox{-0.3ex}{$\m@th\cdot$}}%
                     =}
\makeatother

\makeatletter
\renewcommand{\subsection}{\@startsection{subsection}{2}{\z@}%
  {.5\linespacing\@plus.7\linespacing}%
  {.15\linespacing}%
  {\normalfont\bfseries}}
\makeatother

\begin{document}
\title[]{Categorical Torelli theorem for quintic threefolds}
\author[G.~Ouchi]{Genki Ouchi}

\address{G.~Ouchi, Department of Mathematics, Faculty of Science, Hokkaido University, Kita 10, Nishi 8, Kita-Ku, Sapporo, Hokkaido, 060-0810, Japan}
\email{genki.ouchi@math.sci.hokudai.ac.jp}

\begin{abstract}
We prove that two smooth quintic threefolds over $\bC$ are isomorphic
if their bounded derived categories of coherent sheaves are equivalent.
\end{abstract}

\maketitle{}

\section{Introduction}

By the Bondal--Orlov reconstruction theorem
\cite{bondal-orlov-reconstruction}, two derived-equivalent smooth
projective varieties are isomorphic if the canonical or anticanonical bundle of one of them is ample.
Examples of derived-equivalent but non-isomorphic varieties include
pairs of abelian varieties and pairs of K3 surfaces
\cite{mukai-abelian,mukai-k3,orlov-k3}.
Bridgeland \cite{bridgeland-flops} proved that smooth projective threefolds related by a flop
are derived equivalent. In particular, birational Calabi--Yau threefolds are derived equivalent.
There are also examples of derived-equivalent Calabi--Yau threefolds
that are not birational \cite{caldararu-cy3, borisov-caldararu,hosono-takagi}.
This can happen even when the two threefolds are deformation equivalent
\cite{borisov-caldararu-perry,ottem-rennemo}.
In this paper, we prove the categorical Torelli theorem for quintic threefolds.

\begin{thm}\label{thm:main}
Let $X$ and $Y$ be quintic threefolds. 
If $\Db(X) \simeq \Db(Y)$, then $X \simeq Y$.
\end{thm}

A derived equivalence between Calabi--Yau threefolds induces an
isomorphism of polarized integral Hodge structures on the torsion-free
parts of their third cohomology groups
\cite[p.~857, footnote~8]{addington-donovan-meachan}.
Ottem and Rennemo proved the following more general statement.

\begin{thm}[{\cite[Proposition~2.1]{ottem-rennemo}}]
\label{thm:derived_hodge_isometry}
Let $X$ and $Y$ be smooth projective varieties of odd
dimension $d$. Suppose that $H^k(X,\bZ)=H^k(Y,\bZ)=0$ for every
odd $k<d$. If $\Db(X) \simeq \Db(Y)$, then there is an isomorphism
\[
    H^d(X,\bZ)_{\mathrm{free}}
    \simeq H^d(Y,\bZ)_{\mathrm{free}}
\]
of polarized integral Hodge structures.
\end{thm}

For a quintic threefold $X$, we have $H^1(X,\bZ)=0$ and
$H^3(X,\bZ)$ is torsion-free. Thus derived-equivalent quintic
threefolds have isomorphic polarized integral Hodge structures on
their third cohomology groups.
Voisin proved the following generic Torelli theorem for quintic
threefolds. Let $M$ be the coarse moduli space of quintic threefolds. 

\begin{thm}[{\cite[Theorem~1]{voisin-quintic}}]
\label{thm:integral_generic_torelli}
There exists a non-empty Zariski open subset $U \subset M$ with the
following property\textup{:} for any quintic threefolds $X$ and $Y$ with
$[X] \in U$, if $H^3(X,\bZ)$ and $H^3(Y,\bZ)$ are isomorphic as
polarized integral Hodge structures, then $X \simeq Y$.
\end{thm}

Theorems~\ref{thm:derived_hodge_isometry}
and~\ref{thm:integral_generic_torelli} together imply
Theorem~\ref{thm:main} when $[X] \in U$.
To remove this assumption, we use the following deformation theorem
of Morimura.

\begin{thm}[{\cite[Theorem~1.1]{morimura}}]
\label{thm:morimura}
Let $X$ and $Y$ be derived-equivalent Calabi--Yau manifolds of
dimension $d \geq 3$. There exist a smooth affine variety $S$ with
a closed point $0 \in S$ and smooth projective families
$\cX \to S$ and $\cY \to S$ with $\cX_0 \simeq X$ and
$\cY_0 \simeq Y$, both versal at $0$, such that $\cX_s$ and $\cY_s$
are derived equivalent for a general closed point $s \in S$.
\end{thm}

We apply Theorem~\ref{thm:morimura} to derived-equivalent quintic
threefolds $X$ and $Y$. After an \'etale base change, the two families
give dominant classifying morphisms $\alpha_X,\alpha_Y \colon S \to M$.
For a general closed point $s \in S$, we have $\alpha_X(s) \in U$,
and the fibres $\cX_s$ and $\cY_s$ are derived equivalent.
Theorems~\ref{thm:derived_hodge_isometry}
and~\ref{thm:integral_generic_torelli} give $\cX_s \simeq \cY_s$.
Thus $\alpha_X(s)=\alpha_Y(s)$ for a general closed point $s \in S$.
Since $M$ is separated, $(\alpha_X,\alpha_Y)^{-1}(\Delta_M)$ is closed
in $S$. It contains a dense set of closed points and hence equals $S$.
In particular, $\alpha_X(0)=\alpha_Y(0)$, so $X \simeq Y$.
See Section~\ref{section:proof} for details.

Lin and Zhang \cite{lin-zhang-fm-hypersurfaces} used a similar
deformation argument to prove categorical Torelli for smooth quartic
hypersurfaces in $\bP^{4k+1}$ for every $k \geq 1$. They deform
an equivalence between Kuznetsov components to reduce the proof
to the case of a general quartic. Together with earlier work of
Lin, Rennemo, and Zhang \cite{lin-rennemo-zhang}, their result gives
the following statement for $3 \leq d < N$ with $d \nmid N$.
If $X,Y \subset \bP^{N-1}$ are smooth hypersurfaces of degree $d$,
then an enhanced Fourier--Mukai equivalence between their Kuznetsov
components implies that $X$ and $Y$ are projectively isomorphic.
For further results on categorical Torelli, see the survey of Pertusi
and Stellari \cite{pertusi-stellari} and the more recent papers
\cite{pirozhkov-hypersurfaces, lin-zhang-serre-algebra,
feyzbakhsh-liu-zhang,dell-jacovskis-rota,
elagin-quartic-del-pezzo,jacovskis-meinsma,lin-rampazzo-zhang,
ma-fano-reconstruction,lin-zhang-fano}.

We conclude the introduction with two questions.
\begin{prob}
\begin{enumerate}[label=\textup{(\arabic*)}]
\item Let $X$ be a quintic threefold and $Y$ a smooth projective
variety with $\Db(X) \simeq \Db(Y)$. Is $Y$ necessarily isomorphic
to a quintic threefold? An affirmative answer, together with
Theorem~\ref{thm:main}, would imply that every Fourier--Mukai partner
of $X$ is isomorphic to $X$.

\item Does Theorem~\ref{thm:integral_generic_torelli} hold for all
quintic threefolds $X$ and $Y$, without the assumption $[X] \in U$?
An affirmative answer would give Theorem~\ref{thm:main} directly
from Theorem~\ref{thm:derived_hodge_isometry}.

\end{enumerate}
\end{prob}

\subsection*{Notation and conventions}
\begin{itemize}
\item
We work over $\bC$. A quintic threefold is a smooth hypersurface
of degree five in $\bP^4$. By a Calabi--Yau manifold we mean a
smooth projective variety $X$ with $\omega_X \simeq \cO_X$ and
$H^i(X,\cO_X)=0$ for $0<i<\dim X$.
\item 
We write $\Db(X)$ for the bounded derived category of coherent
sheaves on $X$. A derived equivalence is a $\bC$-linear exact
equivalence of triangulated categories. 
\item
For a finitely generated
abelian group $A$, we write $A_{\mathrm{free}}$ for the quotient
of $A$ by its torsion subgroup.
\item
Given morphisms of schemes $\cX \to S$ and $u \colon T \to S$,
we write $u^*\cX$ or $\cX_T$ for the fibre product $\cX \times_S T$.
\end{itemize}

\subsection*{Acknowledgements}
The author thanks Franco Rota for inspiring conversations about quintic threefolds during his visit to Hokkaido University.
The author was supported by JSPS KAKENHI Grant Numbers
JP25K00906 and JP25K06932.

\subsection*{AI Declaration}
The author used ChatGPT 6.2 Astra to explore proof
strategies, review the arguments, and revise the manuscript.
The author has checked the proofs in the manuscript.

\section{Proof of the main theorem}\label{section:proof}

Let $X$ be a quintic threefold and $S$ a smooth 
variety with a closed point $0 \in S$. 
Suppose that $\pi \colon \cX \to S$ is a smooth projective
morphism with $\cX_0 \simeq X$.
Fix a closed immersion $i_0 \colon \cX_0 \hookrightarrow \bP^4$.

\begin{lem}\label{lem:relative_closed_immersion}
After replacing $S$ by an \'etale neighbourhood of $0$, the closed immersion
$i_0$ extends to a closed immersion $i \colon \cX \hookrightarrow \bP^4_S$
whose image is a smooth relative quintic hypersurface.
\end{lem}
\begin{proof}
The closed immersion $i_0$ defines a closed point $[i_0]$ of the
relative Hom scheme $\underline{\Hom}_S(\cX,\bP^4_S)$ lying over
$0 \in S$.
Since $H^1(\cX_0,i_0^*T_{\bP^4})=0$, it follows from
\cite[Theorem~I.2.10, Lemma~I.2.12.1, and Proposition~I.2.14.2]{kollar-rational-curves},
applied to the graph of $i_0$, that
$\underline{\Hom}_S(\mathcal X,\bP^4_S)$
is smooth over $S$ at $[i_0]$.
By \cite[Tag~057G]{stacks-project}, there exists a locally
closed subscheme $T$ of $\underline{\Hom}_S(\cX,\bP^4_S)$ containing $[i_0]$ such that the
induced morphism $T \to S$ is \'etale.
By the universal property of the Hom scheme, the inclusion
$T \hookrightarrow \underline{\Hom}_S(\cX,\bP^4_S)$ corresponds to
a $T$-morphism $i \colon \cX_T \to \bP^4_T$ whose fibre at
$[i_0] \in T$ is $i_0$. Since $i_0$ is a closed immersion, it follows from
\cite[Tag~05XA]{stacks-project} that, after replacing $T$
by a connected open neighbourhood of $[i_0]$, the morphism $i$ is
a closed immersion. Its image is a smooth relative hypersurface
since $\cX_T$ is smooth over $T$ of relative dimension three.
By flatness and the connectedness of $T$, the degree of the fibres
is constant and hence equals the degree of the central fibre.
\end{proof}

\begin{rem}\label{rem:common_etale_neighbourhood}
Apply Lemma~\ref{lem:relative_closed_immersion} to two families
$\cX \to S$ and $\cY \to S$. The fibre product of the resulting
\'etale neighbourhoods over $S$ gives a common \'etale neighbourhood
$T$ of $0$. By base change, we obtain closed immersions
$\cX_T \hookrightarrow \bP^4_T$ and $\cY_T \hookrightarrow \bP^4_T$.
\end{rem}

Let $B \subset |\cO_{\bP^4}(5)|$ be the open subset parametrizing
quintic threefolds and
$M := B/\operatorname{PGL}_5$ the coarse moduli space of
quintic threefolds, and let $q \colon B \to M$ be the quotient
morphism.
By Lemma~\ref{lem:relative_closed_immersion}, after replacing $S$
by an \'etale neighbourhood of $0$, we may assume that $i_0$
extends to a closed immersion $i \colon \cX \hookrightarrow \bP^4_S$.
Let $\cQ \subset \bP^4 \times B$ be the universal quintic
hypersurface. The relative quintic hypersurface $\cX \subset \bP^4_S$
determines a morphism $\alpha' \colon S \to B$ such that
$\cX = (\alpha')^*\cQ$ as closed subschemes of $\bP^4_S$.
We define the classifying morphism
$\alpha := q \circ \alpha' \colon S \to M$.

\begin{lem}\label{lem:projective_equivalence}
Let $A$ be a local Artinian $\bC$-algebra with residue field $\bC$,
and let $\cX,\cY \subset \bP^4_A$ be smooth relative quintic
hypersurfaces. Then every isomorphism $f \colon \cX \to \cY$ over $A$
is induced by an element of $\operatorname{PGL}_5(A)$.
\end{lem}

\begin{proof}
Let $\m$ be the maximal ideal of $A$, and denote the closed fibres
of $\cX$ and $\cY$ by $\cX_{\bC}$ and $\cY_{\bC}$, respectively.
We first show that the restriction map
$\Pic(\cX) \to \Pic(\cX_{\bC})$ is injective.

Let $\cL$ be a line bundle on $\cX$ whose restriction to
$\cX_{\bC}$ is trivial. Choose a trivialization $\cO_{\cX_{\bC}} \simto \cL|_{\cX_{\bC}}$. Since $H^1(\cX_{\bC},\cO_{\cX_{\bC}})=0$, the
cohomology and base change theorem
\cite[Theorem~III.12.11]{hartshorne} gives a surjection
\[
 H^0(\cX,\cL)\otimes_A \bC
    \to H^0(\cX_{\bC},\cL|_{\cX_{\bC}}).
\]
As $A/\m=\bC$, the restriction map $H^0(\cX,\cL) \to H^0(\cX_{\bC},\cL|_{\cX_{\bC}})$
is also surjective. 
Choose a section $s \in H^0(\cX,\cL)$ whose restriction
corresponds to $1$ under the chosen trivialization.
Since $\m$ is nilpotent, every $x \in \cX$ lies on the
closed fibre $\cX_{\bC}$. For each $x \in \cX$, the class of $s_x$ generates
$\cL_x/\m\cL_x \simeq (\cL|_{\cX_{\bC}})_x$.
By Nakayama's lemma, $s_x$ generates $\cL_x$. 
Thus the morphism $\cO_{\cX} \to \cL$ sending $1$ to $s$
is an isomorphism on every stalk and hence an isomorphism.
This proves the injectivity.

Let $f \colon \cX \simto \cY$ be an isomorphism over $A$. The isomorphism $f_{\bC}$ preserves the hyperplane class, so
$f_{\bC}^*\cO_{\cY_{\bC}}(1) \simeq \cO_{\cX_{\bC}}(1)$.
By the injectivity proved above,
$f^*\cO_{\cY}(1) \simeq \cO_{\cX}(1)$. Fix such an isomorphism.
Since $H^0(\bP^4_A,\cO_{\bP^4_A}(-4))=
H^1(\bP^4_A,\cO_{\bP^4_A}(-4))=0$, the restriction maps
\[
H^0(\cY,\cO_{\cY}(1))
\xleftarrow{\ \sim\ }
H^0(\bP^4_A,\cO_{\bP^4_A}(1))
\xrightarrow{\ \sim\ }
H^0(\cX,\cO_{\cX}(1))
\]
are isomorphisms. Under these identifications, pullback by $f$
and the chosen isomorphism of line bundles induce an automorphism
of $H^0(\bP^4_A,\cO_{\bP^4_A}(1))$. The corresponding projective
transformation restricts to $f$ on $\cX$.
\end{proof}

\begin{prop}\label{prop:dominance}
If $\pi \colon \cX \to S$ is versal at $0$, then the
classifying morphism $\alpha \colon S \to M$ is dominant.
\end{prop}

\begin{proof}
Set $b := \alpha'(0)$. For each $n \geq 0$, let
$B_n := \Spec(\cO_{B,b}/\m_b^{n+1})$, where $\m_b$ is the maximal
ideal of $\cO_{B,b}$, and let $j_n \colon B_n \hookrightarrow B$
be the natural inclusion. Since $\cQ_b \simeq X$, the family
$j_n^*\cQ$ is a deformation of $X$ over $B_n$. 
By versality of
$\pi$ at $0$, there exists a morphism $u_n \colon B_n \to S$
with $u_n(b)=0$ such that
$j_n^*\cQ \simeq u_n^*\cX$
over $B_n$.
Pulling back $\cX=(\alpha')^*\cQ$ along $u_n$, we obtain
$u_n^*\cX=(\alpha'\circ u_n)^*\cQ$
as closed subschemes of $\bP^4_{B_n}$.
By Lemma~\ref{lem:projective_equivalence}, there exists  $g_n \in \operatorname{PGL}_5(B_n)$
such that $g_n \colon j_n^*\cQ \simto (\alpha'\circ u_n)^*\cQ$.
Since $q$ is invariant under $\operatorname{PGL}_5$, we obtain 
$q \circ j_n=q \circ \alpha' \circ u_n=\alpha \circ u_n$.

Let $Z$ be the closure of $\alpha(S)$ in $M$, endowed with its
reduced induced structure. Since $S$ is reduced, $\alpha$ factors
through $Z$. Thus $j_n$ factors through $q^{-1}(Z)$ for every $n$.
Let $\cI$ be the ideal sheaf of $q^{-1}(Z)$ in $B$. Then
$\cI_b \subset \m_b^{n+1}$ for every $n$, so Krull's intersection
theorem gives
\[
    \cI_b \subset \bigcap_{n\geq 0}\m_b^{n+1}=0.
\]
Therefore, there is an open neighbourhood $V$ of $b$ such that $\cI|_V=0$. 
Thus $V \subset q^{-1}(Z)$. Since $B$ is irreducible, $V$ is
dense in $B$, and the closed subset $q^{-1}(Z)$ must equal $B$.
As $q$ is surjective, it follows that $Z=M$. Hence $\alpha$ is
dominant.
\end{proof}

\begin{proof}[Proof of Theorem~\ref{thm:main}]
Let $X$ and $Y$ be derived-equivalent quintic threefolds.
Choose a smooth affine variety $S$, a closed point $0 \in S$, and
smooth projective families $\cX \to S$ and $\cY \to S$ as in
Theorem~\ref{thm:morimura}, with $\cX_0 \simeq X$ and $\cY_0 \simeq Y$.
Both families are versal at $0$, and their fibres $\cX_s$ and $\cY_s$
are derived equivalent for a general closed point $s \in S$.
Applying Lemma~\ref{lem:relative_closed_immersion} and Remark~\ref{rem:common_etale_neighbourhood} to both families,
we may replace $S$ by a common connected \'etale neighbourhood of $0$
and regard $\cX$ and $\cY$ as smooth relative quintic hypersurfaces in
$\bP^4_S$. Denote their classifying morphisms by
$\alpha_X, \alpha_Y \colon S \to M$.

Let $U \subset M$ be the open subset in
Theorem~\ref{thm:integral_generic_torelli}. Since versality at $0$
is preserved under \'etale base change,
Proposition~\ref{prop:dominance} implies that $\alpha_X$ is dominant.
Thus $\alpha_X^{-1}(U)$ is dense in $S$. For a general closed point $s \in S$,
we have $[\cX_s] \in U$ and $\Db(\cX_s) \simeq \Db(\cY_s)$.
By Theorem~\ref{thm:derived_hodge_isometry}, we obtain an isomorphism
$H^3(\cX_s,\bZ) \simeq H^3(\cY_s,\bZ)$
of polarized integral Hodge structures. By
Theorem~\ref{thm:integral_generic_torelli}, it follows that
$\cX_s \simeq \cY_s$. Equivalently, $\alpha_X(s)=\alpha_Y(s)$.
Since $M$ is separated, the diagonal $\Delta_M$ is closed in $M \times M$.
It follows that $(\alpha_X,\alpha_Y)^{-1}(\Delta_M)$ is closed in $S$.
This closed subset contains a dense set of closed points of $S$,
so it equals $S$. In particular, $\alpha_X(0)=\alpha_Y(0)$,
and hence $X \simeq Y$.
\end{proof}

\end{document}